\documentclass[reqno]{amsart}
\usepackage{amsmath,amssymb}
\usepackage{booktabs,mathrsfs,enumerate}
\usepackage{hyperref,framed}
\usepackage[hyperpageref]{backref}
\hypersetup{pdftitle={Alternating quiver Hecke algebras},
  pdfauthor={Clinton Boys and Andrew Mathas},
  colorlinks = true,
  linkcolor =blue,
  anchorcolor = red,
  citecolor = blue,
  urlcolor = blue
}

\renewcommand*{\backref}[1]{}
\renewcommand*{\backrefalt}[4]{%
  \ifcase #1 No citations.% shouldn't happen...
  \or [Page #2.]%
  \else [Pages #2.]%
  \fi%
}

\usepackage{todonotes}

\usepackage{aliascnt}
\def\NewTheorem#1{%
  \newaliascnt{#1}{equation}
  \newtheorem{#1}[#1]{#1}
  \aliascntresetthe{#1}
  \expandafter\def\csname #1autorefname\endcsname{#1}
}
\def\equationautorefname~#1\null{(#1)\null}

\theoremstyle{plain}
\newcounter{mainTheorem}

\newcounter{mainCorollary}
\numberwithin{mainCorollary}{mainTheorem}

\newtheorem{Corollary*}[mainCorollary]{Corollary}

\swapnumbers
\numberwithin{equation}{section}
\NewTheorem{Assumption}
\NewTheorem{Proposition}
\NewTheorem{Theorem}
\NewTheorem{Corollary}
\NewTheorem{Conjecture}
\NewTheorem{Lemma}
\theoremstyle{definition}
\NewTheorem{Definition}
\theoremstyle{remark}
\NewTheorem{Remark}
\NewTheorem{Remarks}

\newaliascnt{Example}{equation}
\aliascntresetthe{Example}

\newenvironment{Example}%
  {\refstepcounter{Example}\trivlist
   \item[\hskip\labelsep\theequation.~\textbf{Example}\space]
   \ignorespaces
  }{\unskip\nobreak\hfil%
    \penalty50\hskip2em\hbox{}\nobreak\hfil$\Diamond$%
    \parfillskip=0pt\finalhyphendemerits=0\penalty-100\endtrivlist
}
\newaliascnt{Examples}{equation}
\aliascntresetthe{Examples}

  {\refstepcounter{Examples}\trivlist
   \item[\hskip\labelsep\theequation.~\textbf{Examples}\space]
   \ignorespaces #1\enumerate
   }{\endenumerate\unskip\nobreak\hfil%
    \penalty50\hskip2em\hbox{}\nobreak\hfil$\Diamond$%
    \parfillskip=0pt\finalhyphendemerits=0\penalty-100\endtrivlist
}

{\catcode`\|=\active
  \gdef\set#1{\mathinner{\lbrace\,{\mathcode`\|"8000%
                                   \let|\midvert #1}\,\rbrace}}
}
\def\midvert{\egroup\,\mid\,\bgroup}

\def\({\big(}
\def\){\big)}

\let\<=\langle
\let\>=\rangle

\usepackage{tikz}
\usetikzlibrary{matrix}

\newcount\tableauRow\newcount\tableauCol
\newcommand\Tableau[2][-2]{
  \begin{tikzpicture}[scale=0.3,draw/.append style={thick,black},baseline=#1mm]
    \tableauRow=0
    \foreach \Row in {#2} {
       \tableauCol=1
       \foreach\k in \Row {
          \draw(\the\tableauCol,\the\tableauRow)+(-.5,-.5)rectangle++(.5,.5);
          \draw(\the\tableauCol,\the\tableauRow)node{\k};
          \global\advance\tableauCol by 1
       }
       \global\advance\tableauRow by -1
    }
  \end{tikzpicture}
}

\def\bi{\mathbf{i}}

\renewcommand{\O}{\mathcal{O}}

\newcommand\RootMr[1][\pm]{\sqrt{\vrule height 2mm width 0pt\smash{M_r'}}}
\newcommand\RootLr[1][\pm]{\sqrt{\vrule height 2mm width 0pt\smash{L_r'}}}

\newcommand{\Sn}[1][n]{\mathfrak{S}_{#1}}

\newcommand{\N}{\mathbb{N}}

\newcommand{\Z}{\mathbb{Z}}

\newcommand{\Zcal}{\mathcal{Z}}

\newcommand{\la}{\lambda}

\DeclareMathOperator{\height}{ht}

\newcommand{\ei}{{\mathbf{i}}}
\newcommand{\ej}{{\mathbf{j}}}

\newcommand{\Y}{{\mathbb{Y}}}

\newcommand{\sgn}{\mathtt{sgn}}

\DeclareMathOperator{\Res}{Res}

\makeindex
\def\Email#1{\email{\href{mailto:#1}{#1}}}
\begin{document}
\title{Alternating quiver Hecke algebras}
\author{Clinton Boys}
\address{School of Mathematics and Statistics F07, University of
Sydney\newline\hspace*{3mm} NSW 2006, Australia.}
\Email{clinton.boys@sydney.edu.au}
%\author{Andrew Mathas}
%\Email{andrew.mathas@sydney.edu.au}

%\AM{I think that 'Alternating quiver Hecke algebras' is a better title.}

\begin{abstract}
For simply-laced quivers, we consider the fixed-point subalgebra of the quiver Hecke algebra under the homogeneous sign map. This leads to a new family of algebras we call alternating quiver Hecke algebras. We give a basis theorem and a presentation by generators and relations which is strikingly similar to the KLR presentation for quiver Hecke algebras. 
\end{abstract}

\maketitle

% fudge for typesetting \LaTeX commands
%   Usage:    \Macro[e]{fred}  -> \begin{fred}...\end{fred}
%             \Macro[b}{fred   -> \fred{...}
%             \Macro{fred}     -> \fred
\newcommand\Macro[2][-]{%
\ifx#1e$\backslash$begin\{#2\}$\dots\backslash$end\{#2\}
\else $\backslash$#2\ifx#1b$\{\dots\}$\fi
\fi}

\newcommand{\Zo}{\mathcal{Z}}

\section*{Introduction}

% !TEX root =thesis.tex

The study of quiver Hecke algebras is a recent development in representation theory, dually motivated from studying representations of cyclotomic Hecke algebras and Khovanov-Lauda-Rouquier (KLR) algebras. Whereas cyclotomic Hecke algebras, being related to symmetric groups and classical Iwahori-Hecke algebras, have been studied for some time, KLR algebras are relative newcomers, and the two ostensibly different families were linked by groundbreaking work of Brundan and Kleshchev \cite{BK:GradedKL}. Because of this deep and meaningful link, it is becoming customary to refer to KLR algebras as quiver Hecke algebras. In recent years much work has been done to study these algebras in their own right. Mathas \cite{MathasSurvey} gives a particularly erudite survey of such developments. 

Quiver Hecke algebras are associative graded algebras whose presentation by generators and relations, and whose grading, depends on the data of a quiver. Particularly interesting is the case when the quiver is simply laced, as quotients of these algebras are isomorphic to cyclotomic Hecke algebras, endowing these algebras with a $\Z$-grading. In this paper, we define a new family of algebras called alternating quiver Hecke algebras. For quiver Hecke algebras whose quivers are simply laced, it is possible to define a homogeneous involution on the quiver Hecke algebra; studying the fixed-point subalgebra under this involution gives rise to our family of algebras. We discuss these algebras using a version of Clifford theory for associative algebras, and construct a homogeneous basis and a presentation by homogeneous generators and relations which are reminiscent of Khovanov and Lauda \cite{KhovLaud:diagI} and Rouquier's \cite{Rouq:2KM} theorems for quiver Hecke algebras. Cyclotomic quotients of these algebras are studied in \cite{BM2} and \cite{BM}.

This paper is organised as follows. We start by giving the definition of quiver Hecke algebras for arbitrary quivers. Then we discuss the Clifford theory for associative algebras that will give the technical mechanism for most of our proofs; this relies on the construction of the opposite quiver. In Chapter 3 we discuss alternating quiver Hecke algebras, giving a basis theorem for this new family of algebras. Finally, in Chapter 4 we prove the main result of this paper, giving a KLR-style presentation for alternating quiver Hecke algebras. 

{\bf Acknowledgements.} The author was supported by an Australian Postgraduate Award. % and the
%second author by the Australian Research Council. 
The author would like to thank his supervisor Andrew Mathas for his guidance throughout his PhD candidature, as well as his thesis examiners for suggesting some improvements to the arguments in this paper. A slightly modified version of this paper appears as Chapter 5 of the author's PhD thesis \cite{BoysThesis} at the University of Sydney. 

\section{Quiver Hecke algebras and opposite quivers}

In this chapter we define quiver Hecke algebras, and give some basic properties. These algebras were introduced by
Khovanov and Lauda~\cite{KhovLaud:diagI} and Rouquier~\cite{Rouq:2KM}.

%Let $e\in\{0,1,3,4,\ldots\}$ and 
Let $\Gamma$ be a quiver with vertex set $I$. Following Kac~\cite{Kac}, to the simply-laced quiver $\Gamma$ we attach the usual
Lie theoretic data of the positive roots $\set{\alpha_i\mid i\in I}$, the
fundamental weights $\set{\Lambda_i|i\in I}$, the non-degenerate
pairing $(\Lambda_i,\alpha_j)=\delta_{ij}$, for $i,j\in I$, and the
{Cartan matrix} $C=(c_{ij})_{i,j\in I}$ where
\begin{equation}\label{cartanmatrixdef}
c_{ij}=\begin{cases}
        2,&\text{if }i=j,\\
        -1,&\text{if $i\leftarrow j$ or $i\rightarrow j$},\\
        0,&\text{otherwise.}
      \end{cases}
\end{equation}
Let $P^+=\bigoplus_{i\in I}\N\Lambda_i$ and $Q^+_\Gamma=\bigoplus_{i\in I}\N\alpha_i$.
The {height} of $\beta=\sum_i a_i\alpha_i\in Q^+_\Gamma$ is the
non-negative integer $\height\beta=\sum_ia_i\in\N$. Fix $n\ge0$ and let
$Q^+_n=\set{\beta\in Q^+_\Gamma\mid\height\beta=n}$. For $\beta\in Q^+_\Gamma$ let
\begin{equation}\label{Ialphadef}
I^\beta=\set{\bi=(i_1,\dots,i_n)\in I^n|\beta=\alpha_{i_1}+\dots+\alpha_{i_n}}.
\end{equation}
In this paper this data plays only a superficial role in describing the combinatorics of the 
grading on the algebras that we consider.

Let $\Zcal$ be an arbitrary unital associative integral domain. 

\begin{Definition}[Khovanov and Lauda \cite{KhovLaud:diagI}
                   and Rouquier~\cite{Rouq:2KM}]
 \label{D:klrdef}
 Let $\Gamma$ be a simply-laced quiver with vertex set $I$ and suppose that  $\beta\in Q^+_\Gamma$. % and $e\in\{0,3,4,\ldots\}$. 
  The
  {quiver Hecke algebra} $\mathcal{R}_\beta(\Gamma)=\mathcal{R}_\beta^\Lambda(\Gamma,\Zcal)$ is the unital associative $\Zcal$-algebra  with generators
  \begin{equation*}
    \{\psi_1,\ldots,\psi_{n-1}\}\cup \{y_1,\ldots,y_n\}
         \cup \{e_\Gamma(\ei)\mid \ei\in I^\beta\}
  \end{equation*}
  and relations
  {\setlength{\abovedisplayskip}{2pt}
   \setlength{\belowdisplayskip}{1pt}
  \begin{xalignat*}{3}
    e_\Gamma(\ei)e_\Gamma(\ej)&= \delta_{\ei_\Gamma\ej_\Gamma}e_\Gamma(\ei),
        &\textstyle\sum_{\ei\in I^\beta}e_\Gamma(\ei)&= 1,\\%\label{E:top3}
  y_re_\Gamma(\ei)&= e_\Gamma(\ei)y_r,& \psi_re_\Gamma(\ei)&= e(s_r\cdot\ei_\Gamma)\psi_r,& y_ry_s&= y_sy_r,
  %\label{E:next3}
  \end{xalignat*}
  \begin{xalignat*}{2}
    \psi_ry_{r+1}e_\Gamma(\ei)&= (y_r\psi_r+\delta_{i_r i_{r+1}})e_\Gamma(\ei),
  & y_{r+1}\psi_re_\Gamma(\ei)&= (\psi_ry_r+\delta_{i_r i_{r+1}})e_\Gamma(\ei),
  \end{xalignat*}
  \begin{align*}
  \psi_ry_s&= y_s\psi_r,&\text{if }s\neq r,r+1,\\
  \psi_r\psi_s&= \psi_s\psi_r,&\text{if }|r-s|>1,
  \end{align*}
  \begin{align}
  \psi_r^2e_\Gamma(\ei)&= \begin{cases}
      0,&\text{if }i_r=i_{r+1},\\
      (y_r-y_{r+1})e_\Gamma(\ei),&\text{if }i_r\to i_{r+1},\\
      (y_{r+1}-y_r)e_\Gamma(\ei),&\text{if }i_r\leftarrow i_{r+1},\\
      e_\Gamma(\ei),&\text{otherwise},
    \end{cases}\label{E:psirsquared}\\
  \psi_r\psi_{r+1}\psi_re_\Gamma(\ei)&=\begin{cases}
    (\psi_{r+1}\psi_r\psi_{r+1}-1)e_\Gamma(\ei),&\text{if }i_r=i_{r+2}\to i_{r+1},\\
    (\psi_{r+1}\psi_r\psi_{r+1}+1)e_\Gamma(\ei),&\text{if }i_r=i_{r+2}\leftarrow i_{r+1},\\
    \psi_{r+1}\psi_r\psi_{r+1}e_\Gamma(\ei),&\text{otherwise,}
  \end{cases}\label{E:deformedbraid}
  \end{align}
  }%
  for $\ei,\ej\in I^\beta$ and all admissible $r$ and $s$. If $n\ge0$
  then the {quiver Hecke algebra} is the algebra $\mathcal{R}_n=\mathcal{R}_n(\Gamma,\Zcal)=\bigoplus_{\beta\in Q^+_\Gamma}\mathcal{R}_\beta(\Gamma,\Zcal)$.
\end{Definition}

%Henceforth, we fix $\Lambda=\Lambda_0$ so that $\Rbeta=\RBeta$.

%We write $\mathcal{R}_n^\Lambda=\mathcal{R}_n^\Zcal$ when we want to emphasise that $\mathcal{R}_n$ is an
%$\Zcal$-algebra.  
One can check easily that the relations in \autoref{D:klrdef} are homogeneous with respect
to the following $\Z$-valued degree function:
\begin{align}
  \deg e_\Gamma(\ei)&= 0, &\text{for all }\ei\in I^n,\nonumber\\
  \deg y_r&= 2,&\text{for }1\le r\le n,\label{degfunc}\\
  \deg \psi_re_\Gamma(\ei)&= -c_{i_r,i_{r+1}},&\text{for $1\le r<n$ and $\ei\in I^n$}.
  \nonumber
\end{align}
Therefore, $\mathcal{R}_n(\Gamma)$ is a $\Z$-graded algebra.

Using \eqref{Ialphadef}, if we define 
\begin{equation}\label{blockidem}
e_\Gamma(\alpha)=\sum_{\ei\in I^\beta}e_\Gamma(\ei)
\end{equation}
for $\beta\in Q^+_\Gamma$ %then the full quiver Hecke algebra $\mathcal{R}_n:=\mathcal{R}_n(\Gamma,\Zo)$ is defined as 
%\begin{equation}\label{blockdecomp}
%\mathcal{R}_n=\bigoplus_{\alpha\in Q_e}\mathcal{R}_\alpha.
%\end{equation}
then the algebras $\mathcal{R}_\beta=e_\Gamma(\beta) \mathcal{R}_ne_\Gamma(\beta)$ are actually the {\em blocks} of $\mathcal{R}_n$ \cite[Corollary 2.11]{KhovLaud:diagI}.% (the {\em blocks} of an associative algebra are its indecomposable two-sided ideals) 
%and the summand $\mathcal{R}_\alpha=\mathcal{R}_\alpha(\Gamma,\Zo)=e_\alpha \mathcal{R}_n e_\alpha$ is a two-sided ideal of the algebra $\mathcal{R}_n$.% which has block decomposition 

%\CB{Define $\psi_\omega$}

For each permutation $\omega\in\Sn$, fix a reduced expression $\omega=s_{i_1}s_{i_2}\cdots s_{i_n}$ and define
\[
\psi_\omega=\psi_{i_1}\psi_{i_2}\cdots\psi_{i_n}.
\]
Importantly, because the braid relations for the elements $\{\psi_1,\psi_2,\ldots,\psi_{n-1}\}$ are more complicated than the symmetric group braid relations, the elements $\psi_\omega$ {\em do} depend on the choice of reduced expression. For $\mathbf{a}=(a_1,a_2,\ldots,a_n)\in \mathbb{N}^n$, write $y^{\mathbf{a}}\in \mathcal{R}_n$ for the monomial  $y_1^{a_1}y_2^{a_2}\cdots y_n^{a_n}$. The following theorem appears as \cite[Theorem 2.5]{KhovLaud:diagI} and \cite[Theorem 3.7]{Rouq:2KM}. %We translate the Khovanov-Lauda proof into our algebraic language below. 

\begin{Theorem}[Basis theorem for quiver Hecke algebras]\label{qHbasis}
Let $\beta\in Q^+$. Then $\mathcal{R}_\beta(\Gamma)$ is a free $\Zo$-algebra with homogeneous basis
\[
\{\psi_\omega y^{\bf a}e_\Gamma(\ei)\mid \omega\in\Sn,\;{\bf a}\in \N^n,\; \ei\in I^\beta\}. 
\]
\end{Theorem}

In this paper we will be interested in the following construction on quivers.

\begin{Definition}
For a quiver $\Gamma$ with vertex set $I$, let $\Gamma'$ be the {\em opposite quiver} whose vertex set is $I$ and which has edges $j\to i$ whenever there is an edge $i\to j$ in $\Gamma$, for $i,j\in I$. 
\end{Definition}

\begin{Remark}
Note that a quiver may be isomorphic to its own opposite under a bijection $\tau:I\to I$. In this case there will be edges $\tau(i)\to \tau(j)$ in $\Gamma'$ for every edge $i\to j$ in $\Gamma$. This is why we label idempotents $e_\Gamma(\ei)$ with the quiver $\Gamma$.%, as $e_\ei_\Gamma$ means something different to $\ei_{\Gamma'}$.  
\end{Remark}

\begin{Example}\label{typeAex}
Fix an integer $e\in\{0,3,4,5,\ldots\}$ and take $\Gamma$ to be  
the quiver $\Gamma_e$ with vertex set $I=\Z/e\Z$ and edges $i\to i+1$, for $i\in
I$.  That is, $\Gamma_e$ is the infinite quiver of type $A_\infty$ if
$e=0$ and the finite quiver of Dynkin type $A_{e-1}^{(1)}$ if $e\geq 3$.
%We exclude the case when $e=2$ because this
%corresponds to $\xi=-1$, which is excluded from \autoref{D:AltHecke}.
Notice that $\Gamma_e'$ is isomorphic to $\Gamma_e$ under the bijection
\begin{align*}
\tau_e:I&\to I\\
i&\mapsto -i\quad\mbox{mod $e$}.
\end{align*}
The results in this paper in the particular case when $\Gamma=\Gamma_e$ are discussed at length in \cite[Chapter 5]{BoysThesis}. This is how the result for general simply-laced $\Gamma$ may be motivated. Throughout this paper we give examples of how our general results work for this special case.
\end{Example}

\section{Clifford theory for associative algebras}

In order to deduce results about our fixed-point subalgebras from the corresponding results for the full algebra, we will use the language of Clifford theory. Clifford theory was initially developed to study the representations of normal subgroups of finite groups \cite{Clifford}. Here we adapt it to cover associative algebras with a $C_2$-graded Clifford system; details in the finite group case can be found in \cite{CRR}, and a slightly different and more general treatment for associative algebras is given in \cite{RR}. Below we write the two elements of the cyclic group $C_2$ as $\{+,-\}$, where signs multiply according to the usual rules. Again, let $\Zcal$ be an arbitrary associative unital integral domain. 
%, or as $\{0,1\}$ when it is convenient to do so; it will be clear from context which convention we are using. 

%We also note that Mitsuhashi referred to these as $\Z_2$-graded Clifford systems in \cite{Mitsuhashi2}. 
\begin{Definition}\label{clifforddef}
Let $A$ be a $\Zo$-algebra. A {\em $C_2$-graded Clifford system} for $A$ is a family $\{A_s\mid s\in C_2\}$ of two $\Zo$-submodules of $A$ such that
\begin{itemize}
\item[(i)] $A_sA_t=A_{st}$, for $s,t\in C_2$;
\item[(ii)] there exists a distinguished central element $\varepsilon\in A$ such that $\varepsilon^2=1$ and $A_+=\varepsilon A_-$;
\item[(iii)] $A=A_+\oplus A_-$; and
\item[(iv)] $1\in A_+$.
\end{itemize}
$A_+$ is called the {\em even} part of the algebra; $A_-$ is called the {\em odd} part. 
\end{Definition}

Clifford theory allows us to give a neat description of the representation theory of the subalgebra $A_+$ given knowledge of the representation theory of $A$. Specifically, for an $A_+$-module $M$, we may twist the $A_+$-action by $\varepsilon$ to define the module $M^\varepsilon$, which is $M$ as an $\O$-module and where $a\in M$ acts as 
\[
a\cdot m=\varepsilon a \varepsilon m. 
\]
Using Definition \ref{clifforddef}(ii) to realise $C_2$ as $\{1,\varepsilon\}$, the {\em inertia group} of $M$ is 
\[
\mathcal{I}(M)=\{x\in C_2\mid M^x\cong M\}\unlhd C_2. 
\]
The size of the inertia group (either 1 or 2) determines the behaviour of an irreducible representation $N$ restricted from $A$ to $A_+$, which we denote by $\Res^A_{A_+}\!N$, in the following sense. We refer the reader to \cite[pp344--345]{CRM} for a proof of the following result, noting that inverting 2 is necessary to prove the direct sum decomposition in (ii). 

\begin{Proposition}[Clifford's theorem for $C_2$-graded associative algebras]\label{clifford}
Let 2 be invertible in $\Zo$ and let $A$ be an associative $\Zo$-algebra with a $C_2$-graded Clifford system.  Let $N$ be an irreducible $A$-module. Then
\begin{itemize}
\item[(i)] If $\mathcal{I}(\Res^A_{A_+}\!N)=C_2$ then $\Res^A_{A_+}\!N\cong \Bigl(\Res^A_{A_+}\!N\Bigr)^\varepsilon$ is an irreducible $A_+$-module. 
\item[(ii)] If $\mathcal{I}(\Res^A_{A_+}\!N)=1$ then $\Res^A_{A_+}\!N=M_+\oplus M_-$ is the direct sum of two irreducible $A_+$-modules, related under the conjugation map (i.e. $M_+=M_-^\varepsilon$). 
\end{itemize}
Moreover, {\em all} irreducible $A_+$-modules arise in one of these two ways. 
\end{Proposition}
%\begin{proof}
%Let $L$ be an irreducible $A_+$-submodule of $\Res^A_{A_+}\!N$. Then for any $x\in C_2$, $xL$ is also an $A_+$-submodule of $\Res^A_{A_+}\!N$ since $\varepsilon^2=1$; moreover the irreducibility of $L$ implies that $xL$ is also an irreducible $A_+$-module. This proves part (i) and most of (ii); for proofs of the fact that this gives an exhaustive list of irreducible $A_+$-modules and the direct sum in part (ii) (which requires inverting 2), we refer the reader to the argument on pp344--345 of \cite{CRM}. 
%\end{proof}

%\CB{Should need a requirement that 2 is invertible or something similar??}

%\CB{Double check the inertia groups are the right way around}

%\AM{Give a reference (or sketch a proof).}

%\CB{Worth working out the general theorem for presentations of Clifford subalgebras and putting it here? The idea definitely doesn't appear in the literature as far as I'm aware. Also, state this better. }

%\section{The sign automorphism}

\section{Alternating quiver Hecke algebras}

%We have seen in \S3.5 how to define alternating cyclotomic Hecke algebras by taking the subalgebra of fixed-points under the hash map, which is an involution defined on cyclotomic Hecke algebras with symmetric multicharge. Our goal now is to generalise the Brundan-Kleshchev isomorphism framework to these subalgebras in level 1, and the correct objects on the other side of our isomorphism theorem are what we call

Given a result concerning symmetric groups, it is natural to ask what this result implies for its index-2 subgroup, the alternating group. It is this construction and line of enquiry which motivates our definition of alternating quiver Hecke algebras, which are related to quiver Hecke algebras in much the same way. Indeed, as the alternating group algebra is nothing more than the subgroup of the symmetric group algebra of points fixed by the sign map \cite[\S2.1]{James}, our alternating quiver Hecke algebras will be subalgebras of fixed points under an analogous map, this time a homogeneous involution. Recall the definition of the quiver Hecke algebra $\mathcal{R}_n(\Gamma)$ from \S1. We now extend the original definition of Kleshchev, Mathas and Ram \cite{KMR:UniversalSpecht} of the graded sign map. %Recall the definition of $I$ from \eqref{Idef} and, for $\ei\in I^n$, of $-\ei\in I^n$ from \eqref{minusidef}. 

\begin{Definition}\label{gsgndef}
The {\em graded sign map} $\sgn:\mathcal{R}_n(\Gamma)\to\mathcal{R}_n(\Gamma')$ is the map defined on generators as
\[
e_\Gamma(\ei)\mapsto e_{\Gamma'}(\ei),\qquad y_r\mapsto -y_r,\qquad \psi_r\mapsto -\psi_r. 
\]
\end{Definition}

\begin{Proposition}\cite{KMR:UniversalSpecht}
The graded sign map is a well-defined homogeneous algebra homomorphism $\mathcal{R}_n(\Gamma) \to \mathcal{R}_n(\Gamma')$.  
\end{Proposition}
\begin{proof}
That the map is homogeneous is clear from its definition; checking it is a well-defined algebra homomorphism amounts to checking it respects the list of relations in Definition \ref{D:klrdef}. Once the notation is correctly interpreted, remembering that in the assignment $e_\Gamma(\ei)\xrightarrow{\sgn} e_{\Gamma'}(\ei)$, the $\ei$ on the right is a sequences of vertices in the opposite quiver, whose edges are reversed, this is a straightforward exercise which we leave to the reader. 
\end{proof}

%We need to determine how the $\sgn$ map restricts to the blocks of the quiver Hecke algebras. For a given $\alpha\in Q^+_\Gamma$ with $\alpha=\alpha_{j_1}+\alpha_{j_2}+\ldots+\alpha_{j_n}$, we can define
%\[
%\alpha'=\alpha_{{j_1}}+\alpha_{{j_2}}+\ldots+\alpha_{{j_n}}\in Q^+_{\Gamma},
%\]
%interpreting it as a root for the opposite quiver $\Gamma'$. 
%
%\begin{Lemma}\label{blocksgn}
%For $\ei\in I^n$, $\ei\in I^\alpha$ if and only if $\ei\in I^{\alpha'}$. 
%\end{Lemma}
%
%The above lemma means that, on the level of blocks, $\sgn$ maps between $\mathcal{R}_\alpha(\Gamma)$ and $\mathcal{R}_{\alpha'}(\Gamma')$ and so is an involution on the block $\mathcal{R}_\alpha(\Gamma)$ if $\alpha=\alpha'$ and $\Gamma=\Gamma'$, or on the direct sum $\mathcal{R}_\alpha(\Gamma)\oplus \mathcal{R}_{\alpha'}(\Gamma')$ of blocks otherwise. In the cases we are  interested in, we can identify $\Gamma$ with $\Gamma'$ via an obvious bijection (see Example \ref{typeAex}). 
%

%\CB{Remark about $\alpha$'s in both quivers.}

\begin{Definition}
Let $\Gamma$ be a simply-laced quiver and $\Gamma'$ its opposite. We write $\tau:I\to I$ for the map for which $i\to j$ in $\Gamma$ if and only if $\tau(j)\to \tau(i)$ in $\Gamma'$. $\tau$ is called the {\em reversal map}. 
\end{Definition}

We can extend the reversal map $\tau$ to $Q^+_\Gamma$ as follows. Let $\alpha=\alpha_{i_1}+\alpha_{i_2}+\ldots+\alpha_{i_n}\in Q^+_\Gamma$; then define $\tau(\alpha)$ to be the root $\alpha_{\tau(i_1)}+\alpha_{\tau(i_2)}+\ldots+\alpha_{\tau(i_n)}$. The following lemma gives the consequences of our notation. 

\begin{Lemma}\label{alphalemma}
There is an isomorphism of algebras $\mathcal{R}_\alpha(\Gamma)\cong \mathcal{R}_{\tau(\alpha)}(\Gamma')$. Under this isomorphism $e_\Gamma(\ei)$ is identified with $e_{\Gamma'}(\tau(\ei))$. 
\end{Lemma}

Notice that, on the level of blocks, the $\sgn$ map takes $\mathcal{R}_\alpha(\Gamma)$ to $\mathcal{R}_{\alpha}(\Gamma')$ and so is an involution on the direct sum algebra $\mathcal{R}_\alpha(\Gamma)\oplus \mathcal{R}_{\alpha}(\Gamma')$. Where it is convenient, we often use Lemma \ref{alphalemma} to abuse notation without comment in the sequel. We are now ready to define alternating quiver Hecke algebras, which we will study for the remainder of this paper. %The definition of alternating quiver Hecke algebras can be motivated from the fact that the alternating group algebra is the fixed-point subalgebra of the symmetric group algebra under the ungraded sign map . 

%\begin{Definition}
%Let $(\Gamma,\tau)$ be a pair of a quiver and an involution $\tau:I\to I$ with one fixed point, identifying $\Gamma$ with its opposite quiver $\Gamma'$. We call such a pair a {\em reversible quiver}. 
%\end{Definition}
%
%
%
%Notice that for a reversible quiver $\Gamma$, we can interpret the graded sign map as an {\em involution} on $\mathcal{R}_n(\Gamma)$. Notice also that the pair $(\Gamma_e,\tau_e)$ in Example \ref{typeAex} is reversible (since $e\neq 2$). 
%
%\begin{Remark}$\;$
%%Using our notation $\alpha\in Q_e$, w
%
%\begin{itemize}
%\item[(i)] We note that there is a particular root $\alpha$ which gives a pathological sequence. If $i_0$ is the fixed point of the map $\tau$, then
%\[
%\ei_0=\underbrace{(i_0,i_0,\ldots,i_0)}_{\text{$n$ times}}\]
%corresponds to the root $\alpha=n\alpha_0$. In some results below we need to exclude this case.
%\item[(ii)] It would be interesting to classify reversible quivers; the author is unaware of any such study in the literature. 
%\end{itemize}
%\end{Remark}
%
%\begin{Example}
%If $\Gamma=\Gamma_e$, then the involution $\tau_e$ has one fixed point, $i_0=0$. So in this case $\ei_0=(0,0,\ldots,0)$. 
%\end{Example}

\begin{Definition}\label{alphaequiv}
We write $[\alpha]$ for the equivalence class of $\alpha\in Q^+_\Gamma= Q^+_{\Gamma'}$ under the equivalence relation $\sim$ generated by $\alpha\sim \beta$ if $\beta=\tau(\alpha)$. We write $(Q^+_{\Gamma})_\tau$ for the set of all equivalence classes.
\end{Definition}

\begin{Definition}
For a simply-laced quiver $\Gamma$ and $\alpha\in Q_\Gamma^+$, the {\em alternating quiver Hecke algebra} is the fixed-point subalgebra 
\[
\mathcal{R}_\alpha(\Gamma/\Gamma')^\sgn=\Bigl[\mathcal{R}_\alpha(\Gamma)\oplus \mathcal{R}_\alpha(\Gamma')\Bigr]^\sgn
\]
of the quiver Hecke algebra under the graded sign map. The alternating quiver Hecke algebra $\mathcal{R}_n^\sgn$ is the algebra
\begin{equation}\label{altqhdecomp}
\mathcal{R}_n^\sgn =\bigoplus_{[\alpha]\in (Q^+_\gamma)_\tau}\mathcal{R}_\alpha(\Gamma/\Gamma')^\sgn
\end{equation}
Note that this decomposition is independent of our choice of equivalence class representative $\alpha\in[\alpha]$ by Lemma \ref{alphalemma}; by the same lemma we can identify $\mathcal{R}_n^\sgn$ with a subalgebra of $\mathcal{R}_n$. 
\end{Definition}

\begin{Remark}
It is important to note the summands appearing on the right-hand side of \eqref{altqhdecomp} in general are not the {\em blocks} of the alternating quiver Hecke algebras, as they are not guaranteed to be indecomposable (although most of them are). It would be an interesting problem to classify the blocks of alternating quiver Hecke algebras. 
\end{Remark}

We also want to work with equivalence classes of $I^n$ under the involution on residue sequences $\ei$ induced by the reversal map $\tau$.  More precisely, for $\ei,\ej\in I^n$ let $\sim$ be the equivalence relation on $I^n$ generated by $\ei\sim \ej$ if $\ei=\tau(\ej)$, where $\tau(\ei)=(\tau(i_1),\tau(i_2),\ldots,\tau(i_n))$ for $\ei=(i_1,\ldots,i_n)$. From this we obtain a partition of the set $I^n$ into equivalence classes of size 1 or 2, depending on whether or not $\ei=\tau(\ei)$. We denote the equivalence class containing a sequence $\ei$ by $[\ei]$. We denote the set of equivalence classes by $I^n_\tau$ and in each equivalence class we choose a representative $\ei_\tau\in[\ei_\tau]$. \label{eiplusdef} Similarly we write $I^\alpha_\tau$ for the set of $\tau$-equivalence classes in $I^\alpha$, noting that $I^\alpha_\tau = I^{\tau(\alpha)}_\tau$.

Using these equivalence classes and distinguished elements, we now define an element $\varepsilon\in \mathcal{R}_n$ which will be very important in studying alternating quiver Hecke algebras and their cyclotomic quotients. %Since we are identifying $\Gamma$ with its opposite $\Gamma'$ under $\tau$, we no longer need to adorn our idempotents $e(\ei)$ with a quiver (although explicitly, $e_\Gamma(\ei)=e_{\Gamma'}(\tau(\ei))$). 
Define $\varepsilon$ to be the element
\begin{equation}\label{varepdef}
\varepsilon=\sum_{\substack{[\ei]\in I^n_\tau\\ \ei_\tau\in [\ei]}}\Bigl(e_\Gamma(\ei_\tau)-e_{\Gamma}(\tau(\ei_\tau))\Bigr)\in \mathcal{R}_n.
\end{equation}
Notice that by Lemma \ref{alphalemma} we can also identify $\varepsilon$ with an element of the direct sum algebra $\mathcal{R}_n(\Gamma)\oplus \mathcal{R}_n(\Gamma')$. It also has the important properties
\[
\sgn(\varepsilon)=-\varepsilon\quad\mbox{and}\quad \varepsilon^2 = 1. 
\]

%\begin{Remark}
%The stipulation that $\tau$ only has one fixed point is important in the definition of $\varepsilon$. It means that $\varepsilon$, when projected onto the blocks of the alternating quiver Hecke algebra, is only zero on the single block corresponding to the root $n\alpha_0$ (since then $\tau((\ei_0)_\tau)=(\ei_0)_\tau$). 
%\end{Remark}
%

Recall the idempotents $e_\Gamma(\alpha)$ from \eqref{blockidem}, which define the blocks of quiver Hecke algebras. We now define new idempotents %If $\alpha\neq\alpha'$, we define new idempotents 
\begin{equation}\label{altblockidem}
e_{[\alpha]}=e_\Gamma(\alpha) + e_{\Gamma'}(\alpha),
\end{equation}
noting that since $[e_\Gamma(\alpha)]^\sgn=e_{\Gamma'}(\alpha)$, $e_{[\alpha]}$ is $\sgn$-invariant; in particular by Lemma \ref{alphalemma}, $e_{[\alpha]}$ is independent of the choice of equivalence class representative of $[\alpha]$. %Any given equivalence class $[\ei]\in I^n_\tau$ will contain a sequence from $I^\alpha$ and $I^{\alpha'}$; we write $I^{[\alpha]}_\tau$ %$I^{\alpha/\alpha'}_\sim$ 
%for the collection of all such equivalence classes.
% and if $\alpha=\alpha'$ we write $I^\alpha_\sim$. 
Finally, for an equivalence class $[\ei]\in I^n_\tau$ we write 
\begin{equation}\label{esquaredef}
e[\ei]=e_\Gamma(\ei)+e_{\Gamma'}(\ei)
\end{equation} 
which again does not depend on the choice of equivalence class representative. 

%\CB{Example of how this works in the $\Gamma=\Gamma_e$ case goes here.}

%\begin{Example}
%When $\Gamma=\Gamma_e$ for $e\in\{0,3,4,\ldots,\}$, for a residue sequence $\ei$ with $\ei\neq \underbrace{(0,0,\ldots,0)}_{\text{$n$ zeroes}}$, $[\ei]=\{\ei,-\ei\}$ and so $\varepsilon$ is a sum over all residue sequences with half assigned minus signs; for example we might choose the representative $\ei_\tau$ to be the greatest in lexicographic order in its equivalence class. 
%\end{Example}

%\CB{Important: need to mention how the equivalence classes work blockwise, as we need to consider both $I^\alpha$ and $I^{\alpha'}$ in order for this to make sense..}

We now give a basis theorem for alternating quiver Hecke algebras analogous to Theorem \ref{qHbasis} for quiver Hecke algebras. Recall that $(Q^+_{\Gamma})_\tau$ is the set of all equivalence classes $[\alpha]$. For ${\bf a}=(a_1,a_2,\ldots,a_n)\in\mathbb{N}^n$, let $\left|\mathbf{a}\right|=\sum_{i=1}^n a_i$. 

%\CB{Define $e[\ei]$ earlier}

\newcommand\SetBox[3][60]{\Bigg\{\ #2\ \Bigg| \ \vcenter{\hsize#1mm\centering #3}\Bigg\}}
%    \[    \SetBox[80]{C}{$C$ is a removable $j$-strip of length at most $\gamma_j$ such
%                 that $\row\alpha^\bnu_\blam(j)<\row C\leq \row\alpha^\bnu_\bmu(j)$,
%                 for $j\in J^*$}
%  \]

\begin{Theorem}[Basis theorem for alternating quiver Hecke algebras]\label{aqHbasis}
Let $\Gamma$ be a simply-laced quiver and let 2 be invertible in $\Zo$. For $[\alpha]\in (Q^+_{\Gamma})_\tau$ the alternating quiver Hecke algebra $\mathcal{R}_{\alpha}(\Gamma/\Gamma')^\sgn$ has homogeneous basis
\[
\SetBox[75]{\psi_\omega y^{\mathbf{a}}\varepsilon^b e[\ei]}{$\omega\in\Sn,\quad\mathbf{a}\in\N^n,\quad [\ei]\in I^{\alpha}_\tau$, $b\in \{0,1\},\quad \ell(\omega)+\left|\mathbf{a}\right|+b\equiv 0 \; \mbox{mod 2}$}
%\Bigl\{\psi_\omega y^{\mathbf{a}}\varepsilon^be[\ei]&\;\Bigr|\; \omega\in\Sn,\mathbf{a}\in\N^n, [\ei]\in I^{[\alpha]}_\tau, b\in \{0,1\},\\
%&\Bigl.\qquad\qquad\qquad\qquad \ell(\omega)+\left|\mathbf{a}\right|+b\equiv 0 \; \mbox{mod 2}\Bigr\}\\
%\{\psi_\omega y_1^{a_1}\cdots y_n^{a_n}\varepsilon\sum_{\ej\in [\ei]}e(\ej)&\mid\omega\in\Sn,a_1,\ldots,a_n\in\N,[\ei]\in I^{\alpha/\alpha'}_\sim,\\
%&\qquad\qquad\qquad \ell(\omega)+\sum_ia_i\equiv 1\;\mbox{mod 2}\}. 
%\end{align*}
\]
%\CB{Retypeset this nicer}
%\item[(ii)] If $\alpha\in Q^\pm_{\Gamma/\Gamma'}$ then the alternating quiver Hecke algebra $\mathcal{R}_\alpha(e,\O)^\sgn$ has homogeneous basis
%\begin{align*}
%\{\psi_\omega y_1^{a_1}\cdots y_n^{a_n}\sum_{\ej\in [\ei]}e(\ej)&\mid \omega\in\Sn,a_1,\ldots,a_n\in\N, [\ei]\in I^{\alpha}_\sim,\\
%&\qquad\qquad\qquad \ell(\omega)+\sum_i a_i\equiv 0 \; \mbox{mod 2}\}\;\cup\\
%\{\psi_\omega y_1^{a_1}\cdots y_n^{a_n}\varepsilon\sum_{\ej\in [\ei]}e(\ej)&\mid\omega\in\Sn,a_1,\ldots,a_n\in\N,[\ei]\in I^{\alpha}_\sim,\\
%&\qquad\qquad\qquad \ell(\omega)+\sum_ia_i\equiv 1\;\mbox{mod 2}\}. 
%\end{align*}
%\item[(ii)] If $\alpha=n\alpha_0$ then the alternating quiver Hecke algebra $\mathcal{R}_{n\alpha_0}^\sgn$ has homogeneous basis 
%\[
%\{\psi_\omega y^{\mathbf{a}}e\mid \omega\in\Sn, \;\;\mathbf{a}\in \N^n, \;\;\ell(\omega)+\left|\mathbf{a}\right|\equiv 0\;\mbox{mod\; 2}\},
%\]
%where $e=e(\ei_0)$.%\underbrace{(0,0,\ldots,0)}_{\text{$n$ zeroes}}$. 
%\end{itemize}
%provided $2$ is invertible in $\O$. 
\end{Theorem}

%\CB{Possible to state this simultaneously for $Q^\pm_{\Gamma/\Gamma'}$?}

\begin{proof}
Let $[\alpha]\in (Q^+_{\Gamma})_\tau$. To observe that the specified set spans $\mathcal{R}_{\alpha}(\Gamma/\Gamma')^\sgn$, let us write an arbitrary element of $\mathcal{R}_\alpha(\Gamma)\oplus \mathcal{R}_{\alpha}(\Gamma')\cong \mathcal{R}_\alpha(\Gamma)\oplus\mathcal{R}_{\tau(\alpha)}(\Gamma)$ as a finite sum 
\[
x=\sum_{\omega,\mathbf{a},\ei}\la_{\omega,\mathbf{a},\ei}\psi_\omega y^{\mathbf{a}}e_\Gamma(\ei),
\]
for $\omega\in\Sn$, $\mathbf{a}\in \N^n$, $\ei\in I^\alpha\cup I^{\tau(\alpha)}$ and $\lambda_{\omega,\mathbf{a},\ei}\in \Zo$ using Theorem \ref{qHbasis} (note we are implicitly using Lemma \ref{alphalemma} here to write this as an element of $\mathcal{R}_n(\Gamma)$). In order that $x$ be an element of $\mathcal{R}_{\alpha}^\sgn$ we require $x^\sgn=x$, i.e. that
\[
\sum_{\omega,\mathbf{a},\ei}(-1)^{\ell(\omega)+\left|\mathbf{a}\right|}\la_{\omega,\mathbf{a},\tau(\ei)}\psi_\omega y^{\mathbf{a}}e_\Gamma(\ei)=\sum_{\omega,\mathbf{a},\ei}\la_{\omega,\mathbf{a},\ei}\psi_\omega y^{\mathbf{a}}e_\Gamma(\ei)
\]
which implies that the given elements span by applying Theorem \ref{qHbasis} and equating coefficients. For linear independence, take a linear combination
\[
\sum_{[\ei]\in I^{\alpha}_\tau}\Bigl(\sum_{\substack{\omega\in\Sn,\mathbf{a}\in\N^n\\ \ell(\omega)+\left|\mathbf{a}\right|\equiv 0 \; \text{mod 2}}}\lambda_{\mathbf{a},\omega,\ei}\psi_\omega y^{\mathbf{a}}e[\ei]+\sum_{\substack{\omega\in\Sn,\mathbf{a}\in\N^n\\ \ell(\omega)+\left|\mathbf{a}\right|\equiv 1 \; \text{mod 2}}}\lambda_{\mathbf{a},\omega,\ei}\psi_\omega y^{\mathbf{a}}\varepsilon e[\ei]\Bigr)=0
\]
for coefficients $\la_{\mathbf{a},\omega,\ei}\in \Zo$ and project the above sum separately onto the idempotents $e_\Gamma(\ei)$ and $e_{\Gamma'}(\ei)$ for each $\ei\in I^\alpha$; %then since $\ei\neq \tau(\ei)$ for any idempotents in this case, 
we thus obtain sums of basis vectors of the form in Theorem \ref{qHbasis} which are linearly independent. 
%The proof for $\alpha\in Q^\pm_{\Gamma/\Gamma'}$ is identical to the above; since $\ei\neq -\ei$ for any $\ei\in I^\alpha$ with $\alpha\in Q_e^{\pm}$, the linear independence argument above still holds. 
%For $\alpha=n\alpha_0$, the same argument as above shows the set spans the $\sgn$-invariant subalgebra. Since this basis is a subset of the basis from Theorem \ref{qHbasis}, linear independence, and hence the result, follows. 
\end{proof}

%\CB{Example of what this looks like in the $\Gamma=\Gamma_e$ case goes here.}

%\CB{Change notation above to equivalence classes}
%\CB{Will need a different proof for the case $\alpha=n\alpha_0$ and I think it does still work. }

%For the remaining summand, the homogeneous basis is slightly simpler. 
%
%\begin{proposition}\label{basispath}
%\end{proposition}
%\begin{proof}
%The same argument as above shows the set spans the $\sgn$-invariant subalgebra. Since this basis is a subset of the basis from Theorem \ref{qHbasis}, linear independence and the result follow. 
%\end{proof}

%Under some slightly more restrictive conditions on the root $\alpha$ corresponding to the block, 

We can obtain a $C_2$-graded Clifford decomposition for direct sums of blocks of quiver Hecke algebras using the element $\varepsilon$ defined in \eqref{varepdef}.%, provided $\alpha$ is not equal to the ``pathological'' root $n\alpha_0$.  

\begin{Proposition}\label{clifforddecomp}
Let $\Gamma$ be a simply-laced quiver and let $2$ be invertible in $\Zo$.% and suppose that $\alpha \neq n\alpha_0$. 
%\begin{itemize}
%\item[(i)] Let $\alpha\in Q^+_{\Gamma/\Gamma'}$.
 Then we have the $C_2$-graded Clifford decomposition
\[
\bigoplus_{\beta\in[\alpha]} \mathcal{R}_\beta \cong \mathcal{R}_{\alpha}(\Gamma/\Gamma')^\sgn \oplus \varepsilon \mathcal{R}_{\alpha}(\Gamma/\Gamma')^\sgn. 
\]
%\mathcal{R}_\alpha\oplus \mathcal{R}_{\alpha'}\cong \mathcal{R}_{\alpha/\alpha'}^\sgn\oplus \varepsilon\mathcal{R}_{\alpha/\alpha'}^\sgn. 
%\]
%\item[(ii)] Let $\alpha\in Q^\pm_{\Gamma/\Gamma'}$. Then we have the $C_2$-graded Clifford decomposition
%\[
%\mathcal{R}_\alpha \cong \mathcal{R}_\alpha^\sgn\oplus \varepsilon \mathcal{R}_\alpha^\sgn. 
%\]
%\end{itemize}
\end{Proposition}

%Moreover, we have the $C_2$-graded Clifford decomposition
%\[
%\mathcal{R}_n\cong \mathcal{R}_n^\sgn \oplus \varepsilon\mathcal{R}_n^\sgn.
%\]\\
\begin{proof}
To demonstrate the Clifford decomposition, we check the requirements in Definition \ref{clifforddef}. Condition (i) follows easily; since $\varepsilon^\sgn=-\varepsilon$, if $x\in \mathcal{R}_n^\sgn$ and $y=\varepsilon z\in \varepsilon \mathcal{R}_n^\sgn$ then $(xy)^\sgn=-xy$ giving the required multiplicative property. Conditions (ii) and (iv) follow by definition, so it remains to demonstrate the direct sum decomposition. %, which follows from the same argument as in Theorem \ref{ungradedaltdim}: 
Since any $x\in\mathcal{R}_\alpha\oplus\mathcal{R}_{\alpha'}$ may be written as $x=\sum_{\ei\in I^\alpha\cup I^{\alpha'}}xe_\Gamma(\ei)=\sum_{\ei\in I^\alpha}x\Bigl[e_\Gamma(\ei)+e_\Gamma(\tau(\ei))\Bigr]$, and since $2$ is invertible in $\Zo$, we can write
\[
x=\frac12\sum_{\ei\in I^\alpha}(x+x^\sgn)\Bigl[e_\Gamma(\ei)+e_\Gamma(\tau(\ei))\Bigr]+ \frac12\varepsilon\sum_{\ei\in I^\alpha}(x-x^\sgn)\Bigl[e_\Gamma(\ei)+e_\Gamma(\tau(\ei))\Bigr]
\]
which gives the required decomposition.% provided $\frac12\in\Zo$. %The argument for $\alpha\in Q^\pm_{\Gamma/\Gamma'}$ is identical. 
\end{proof}

%\CB{Check this is indeed the only weight $\alpha$ for which $\varepsilon e_\alpha=0$ -- if not I guess we can just write ``assume $\varepsilon e_\alpha\neq 0$''...}
%
%\CB{If $\alpha=\alpha'$, exactly the same argument will give that $\mathcal{R}_\alpha=\mathcal{R}_\alpha^\sgn\oplus\varepsilon \mathcal{R}_\alpha^\sgn$}

%\CB{Not quite finished with the Clifford decomposition.}

\begin{Remark}\label{clifffailrk}$\;$
%\begin{itemize}
%\item[(i)] We exclude the case $\alpha=n\alpha_0$ in Proposition \ref{clifforddecomp} because in that case the projection of $\varepsilon$ onto the block $\mathcal{R}_\alpha$ is zero, since if $\alpha=n\alpha_0$, there is a single idempotent $e_{n\alpha_0}=e(0,0,\ldots,0)$ which maps to itself under $\sgn$ so $\varepsilon e_{n\alpha_0}=0$.% and condition (ii) fails.
%\item[(i)] For any $i\in I$, the algebra $\mathcal{R}_{n\alpha_i}$ is called a {\em nil-Hecke algebra}. Our algebra $\mathcal{R}_{n\alpha_0}^\sgn$, which behaves slightly differently to the other summands of the alternating quiver Hecke algebra from \eqref{aqHasummands}, may be thought of as an alternating nil-Hecke algebra. %More details on nil-Hecke algebras in the framework we use here can be found in \cite[\S2.5]{MathasSurvey}.
%\item[(ii)] 
The reader may have noticed the dependence of our element $\varepsilon$ on choices of equivalence class representatives for each class $[\ei]\in I^n_\tau$. Indeed, there are many choices of Clifford element which give rise to different Clifford decompositions, much in the way that different choices of coset representatives give rise to different Clifford decompositions of representations of finite groups with respect to normal subgroups \cite{CRM}. 
%\end{itemize}
\end{Remark}

\begin{Example}
We finish this chapter with an example of our results when $\Gamma=\Gamma_e$. First suppose $n=1$ and that $2<e<\infty$. Set $I=\Z/e\Z$ and notice that in this case our quiver Hecke algebra $\mathcal{R}_1(\Gamma_e,\Zcal)$ has a particularly simple presentation:
\begin{align*}
\mathcal{R}_1(\Gamma_e)&=\Bigl\langle y,e(i)\;\Bigr|\;\Bigl. i\in I,\; ye(i)=e(i)y, \;\sum_{i\in I}e(i)=1, \;e(i)e(j)=\delta_{ij}e(i)\Bigr\rangle\\
&\cong \bigoplus_{i\in I}\Zo[y]e(i). 
\end{align*}
where $\deg y =2$ and $\deg e(i)=0$. This has the obvious homogeneous basis
\[
\{y^ke(i)\mid k\geq 0, i\in I\}. 
\]

Suppose now that $e=3$ and that $\mathrm{char}(\Zo)\neq 2$; then it is not too hard to see that the following is a homogeneous basis for the fixed-point subalgebra $\mathcal{R}_1(\Gamma_e)^\sgn$, where we have grouped basis vectors by degree:
\begin{alignat*}{2}
\deg 0:&\quad 2e(0), \;\;e(1)+e(2)\\
\deg 2:&\quad y\Bigl(e(1)-e(2)\Bigr)\\
\deg 4:&\quad 2y^2 e(0), \;\;y^2\Bigl(e(1)+e(2)\Bigr)\\
\deg 6:&\quad y^3\Bigl(e(1)-e(2)\Bigr)\\
&\quad\qquad \vdots
\end{alignat*}
where $y=y_1$. For general $e>2$, we have the homogeneous basis 
\[
\left\{y^k \Bigl(e(i)+e(-i)\Bigr)\mid k \;\mbox{even}\right\}\cup\left\{y^k\Bigl(e(i)-e(-i)\Bigr)\mid k\;\mbox{odd}, \; i\neq 0\right\}
\]
for $\mathcal{R}_n(\Gamma_e)^\sgn$. Note that
\begin{align*}
\dim \Bigl([\mathcal{R}_n(\Gamma_e)^\sgn]_{2k}&\oplus [\mathcal{R}_n(\Gamma_e)^\sgn]_{2(k+1)}\Bigr)\\
&\qquad =\frac12\dim \Bigl([\mathcal{R}_n(\Gamma_e)]_{2k}\oplus [\mathcal{R}_n(\Gamma_e)]_{2(k+1)}\Bigr)
\end{align*}
for all $k>0$. 
%%Finally the reader should check that $\mathcal{R}_1(\Gamma_e)e^+([1])$ has the $C_2$-Clifford decomposition 
%%\begin{align*}
%%\mathcal{R}_1(3,\Zo)e(1)&\oplus \mathcal{R}_1(3,\Zo)e(2)\cong \Bigl[\mathcal{R}_1(3,\Zo)e(1)\oplus \mathcal{R}_1(3,\Zo)e(2)\Bigr]^\sgn \\
%%&\qquad\qquad\qquad\qquad\quad\oplus \varepsilon\Bigl[ \mathcal{R}_1(3,\Zo)e(1)\oplus \mathcal{R}_1(3,\Zo)e(2)\Bigr]^\sgn,
%%\end{align*}
%%where $\varepsilon=e(1)-e(2)$, and that the transition matrix between this basis and the basis from Example \ref{affine1} requires inverting 2.
\end{Example}

%\begin{Example}
%It is worth giving an example of how a Clifford decomposition like that in Proposition \ref{clifforddecomp} is impossible when $\alpha=n\alpha_0$, as mentioned in Remark \ref{clifffailrk}(i). Indeed, let $n=1$ and consider the block $\mathcal{R}_0=e(0)\mathcal{R}_1$ of the algebra from Example \ref{affine1}. Then $\mathcal{R}_0=\Zo[y]e(0)$ so $\mathcal{R}_0^\sgn=\Zo[y^2]e(0)$ since $e(0)^\sgn=e(0)$. Since $\varepsilon e(0)=0$ we do not have a Clifford decomposition in this case. We can also easily compute the basis $\{1,y^2,y^4,\ldots\}$ from Theorem \ref{aqHbasis}(iii). 
%\end{Example}

\section{Generators and relations for alternating quiver Hecke algebras}

Using the $C_2$-graded Clifford decomposition from Theorem \ref{aqHbasis}, we can give a presentation for the alternating quiver Hecke algebra $\mathcal{R}_n(\Gamma/\Gamma')^\sgn$ by homogeneous generators and relations, akin to the Khovanov-Lauda presentation for the quiver Hecke algebras (Definition \ref{D:klrdef}). The following elements will play the role of the generators from Definition \ref{D:klrdef} for alternating quiver Hecke algebras.% in the alternating case. 

\begin{Definition}[Generators for alternating quiver Hecke algebras]\label{aqHagens}
Let $\Gamma$ be a simply-laced quiver and $\mathcal{R}_n(\Gamma/\Gamma')^\sgn$ be an alternating quiver Hecke algebras. For $1\leq r<n$, define $\Psi_r=\psi_r\varepsilon\in \mathcal{R}_n(\Gamma/\Gamma')^\sgn$, and for $1\leq r\leq n$, define $\mathcal{Y}_r=y_r\varepsilon\in\mathcal{R}_n(\Gamma/\Gamma')^\sgn$. Finally, for $[\ei]\in I^n_\tau$ recall the definition of $e[\ei]\in\mathcal{R}_n(\Gamma/\Gamma')^\sgn$ from \eqref{esquaredef}. 
%for $[\ei]\in I^n_\tau$ define 
%\[
%e[\ei]= \sum_{\ej\in [\ei]}e(\ej).
%\]
%We will also need the notation $e^-[\ei]=\varepsilon e[\ei]$ for $[\ei]\in I^n_\tau$. 
%Note that $e^-(\ei)=\varepsilon e^+(\ei)$ for $\ei\in I^n_+$. 
\end{Definition}
%\begin{remark}
%We also allow for $e^+(\ei)$ and $e^-(\ei)$ for $\ei\in I^n_-$, however we never need to use them due to the following simple relations obtained from \eqref{varepdef}:% mean we never have to use these though:
%\begin{align}
%e^+(\ei)&= e^+(-\ei)\label{eminus1}\\
%e^-(\ei)&= -e^-(-\ei)\label{eminus2}
%\end{align}
%for all $\ei\in I^n_-$. 
%\end{remark}
\newcommand{\varep}{\varepsilon}

%\CB{Is it a problem that the {\em generators} depend on our choices of equivalence class representatives?}

%\begin{remark}
%The reader will note that the generators we have defined above depend on the choice of representatives $\ei_\tau\in [\ei]$ that we chose to define the Clifford element $\varep$. We will see as a consequence of Corollary \ref{Z2decomp} that the algebra is actually independent of these choices. 
%\end{remark}

As a corollary to Theorem \ref{aqHbasis} we see that the elements defined above generate the alternating quiver Hecke algebras. 

\begin{Corollary}\label{aqHagenscor}
Let $\Gamma$ be a simply-laced quiver, let $n>1$ and $e>2$ and suppose 2 is invertible in $\Zo$. For $[\alpha]\in (Q_\Gamma^+)_\tau$, the alternating quiver Hecke algebra $\mathcal{R}_{\alpha}(\Gamma/\Gamma')^\sgn$ is generated by the collection of elements
\[
\{\Psi_1,\Psi_2,\ldots,\Psi_{n-1}\}\cup\{\mathcal{Y}_1,\mathcal{Y}_2,\ldots,\mathcal{Y}_n\}\cup\{e[\ei]\mid [\ei]\in I^{[\alpha]}_\tau\}. 
\]
%\item[(ii)] Let $\alpha\in Q^\pm_{\Gamma/\Gamma'}$. The alternating quiver Hecke algebra $\mathcal{R}_\alpha^\sgn$ is generated by the collection
%\[
%\{\Psi_1,\Psi_2,\ldots,\Psi_{n-1}\}\cup\{\mathcal{Y}_1,\mathcal{Y}_2,\ldots,\mathcal{Y}_n\}\cup\{e^+([\ei])\mid [\ei]\in I^\alpha_\sim\}. 
%\]
%\item[(ii)] If $\alpha=n\alpha_0$, the alternating quiver Hecke algebra $\mathcal{R}_{n\alpha_0}^\sgn$ is generated by all even products of the generators $\{\psi_1,\psi_2,\ldots,\psi_{n-1}\}\cup\{y_1,y_2,\ldots,y_n\}$ and $\{e\}$ where $e=e(\ei_0)$.  
%\end{itemize}
\end{Corollary}
\begin{proof}
Suppose $[\alpha]\in (Q_\Gamma^+)_\tau$. We proceed to write each basis vector from Theorem \ref{aqHbasis} in terms of the proposed generators; since the $y_r$ commute with the $e_\Gamma(\ei)$ for all admissible $r$ and $\ei$ by the relations in Definition \ref{D:klrdef}, and since $\varep^2=1$, 
\[
y_1^{a_1}\cdots y_n^{a_n}=\left\{\begin{array}{ll} \mathcal{Y}_1^{a_1}\cdots \mathcal{Y}_n^{a_n},&\mbox{if $\left|\mathbf{a}\right|\equiv 0\;\mbox{mod}\;2$}\\ \varep\mathcal{Y}_1^{a_1}\cdots \mathcal{Y}_n^{a_n},&\mbox{if $\left|\mathbf{a}\right|\equiv 1\;\mbox{mod}\;2$}\end{array}\right.
\]
Moreover, since $\psi_\omega e_\Gamma(\ei)=e_\Gamma(\omega\cdot\ei)\psi_\omega$ by the relations in Definition \ref{D:klrdef}, 
\[
\psi_\omega= \left\{\begin{array}{ll} \Psi_\omega,&\mbox{if $\ell(\omega)\equiv 0\;\mbox{mod}\;2$}\\ \Psi_\omega\varep,&\mbox{if $\ell(\omega)\equiv 1\;\mbox{mod}\;2$}\end{array}\right.
\]
There are four cases to consider; we just give two as illustration of the method of proof. If $\left|\mathbf{a}\right|\equiv 0\;\mbox{mod}\;2$ and $\ell(\omega)\equiv 0\;\mbox{mod}\;2$, then the basis vector $\psi_\omega y_1^{a_1}\cdots y_n^{a_n}e[\ei]$ is equal to $\Psi_{i_1}\cdots\Psi_{i_d}\mathcal{Y}_1^{a_1}\cdots\mathcal{Y}_n^{a_n} e[\ei]$ where $\omega=s_{i_1}\cdots s_{i_d}$. If $\left|\mathbf{a}\right|\equiv 1\;\mbox{mod}\;2$ and $\ell(\omega)\equiv 0\;\mbox{mod}\;2$ then the basis vector 
\begin{align*}
\psi_\omega y_1^{a_1}\cdots y_n^{a_n}[e_\Gamma(\ei)-e_\Gamma(\tau(\ei))]&= \psi_\omega y_1^{a_1}\cdots y_n^{a_n}\varep e[\ei]\\
&= \Psi_{i_1}\cdots\Psi_{i_d} \varep \mathcal{Y}_1^{a_1}\cdots \mathcal{Y}_n^{a_n}\varep e[\ei]\\
&= \Psi_{i_1}\cdots\Psi_{i_d} \mathcal{Y}_1^{a_1}\cdots \mathcal{Y}_n^{a_n} e[\ei]
\end{align*}
since $\varep$ commutes with the $\mathcal{Y}_r$ and squares to 1. %The other two cases are similar. 
%This argument shows the given elements generate the summands $\mathcal{R}_{\alpha/\alpha'}^\sgn$ for $\alpha\in Q^+_{\Gamma/\Gamma'}$ and $\mathcal{R}_\alpha^\sgn$ for $\alpha\in Q^\pm_{\Gamma/\Gamma'}$. 
%When $\alpha=n\alpha_0$, the argument is essentially the same, and easier than, the case above and we leave this to the reader.  % and so the generators $\Psi_1,\Psi_2,\ldots,\Psi_{n-1},\mathcal{Y}_1,\mathcal{Y}_2,\ldots,\mathcal{Y}_n$ are all zero in $\mathcal{R}_{n\alpha_0}^\sgn$. However it is clear by inspecting the generators and relations for the block $\mathcal{R}_{n\alpha_0}$ of the full quiver Hecke algebra that the specified collection of elements in part (iii) generate the fixed-point subalgebra. 
\end{proof}

%\CB{Change the above to equivalence classes}
%\CB{Blockwise?}
%\CB{Write down the result and proof separately for the different $\alpha$ cases? Include more details of the $n\alpha_0$ proof?}

Our goal now is to obtain a set of {\em relations} for the generators from Corollary \ref{aqHagenscor}. Let us start by defining a new abstract algebra with a slightly different presentation and the additional structure of a $\Z\times C_2$-grading; as in \S5.1, let us realise $C_2$ as the sign group $\{+,-\}$ with usual multiplication of signs.%: $+-=-+=-$, $++=--=+$. 

%\AM{Be consistent with choosing $C_2$ or $C_2$}

\begin{Definition}\label{twisteddef}
Let $n\geq 0$. The {\em signed quiver Hecke algebra} of type $\Gamma$ and corresponding to $[\alpha]\in (Q^+_{\Gamma})_\tau$ is the unital associative $\Zo$-algebra $\mathcal{R}_{[\alpha]}(\Gamma)'=\mathcal{R}_{[\alpha]}'(\Gamma,\Zcal)$ with generators
\[
\{\psi_1',\psi_2',\ldots,\psi_{n-1}'\}\cup\{y_1',y_2',\ldots,y_n'\}\cup\Bigl\{\varepsilon_a(\ei)\mid \ei\in \bigcup_{\beta\in [\alpha]}I^\beta,\; a\in C_2\Bigr\}
\]
subject to the relations
\begin{align*}
\varep_a(\ei)\varep_b(\ej)&= \delta_{\ei\ej}\varep_{ab}(\ei),\qquad \sum_{\ei\in I^{[\alpha]}}\varep_+(\ei)= 1\\
\varep_a(\ei)&= a\varep_a(\tau(\ei))\\
y_r'\varep_a(\ei)&= \varep_a(\ei)y_r'\\
\psi_r'\varep_a(\ei)&= \varep_a(s_r\cdot\ei)\psi_r'\\
y_r'y_s'\varep_-(\ei)&= y_s'y_r'\varep_-(\ei)\\
\psi_r'y_s'\varep_-(\ei)&= y_s'\psi_r'\varep_-(\ei)\quad\mbox{if $s\neq r,r+1$}\\
\psi_r'\psi_s'\varep_-(\ei)&= \psi_s'\psi_r'\varep_-(\ei)\quad\mbox{if $\left|r-s\right|>1$}\\
\psi_r'y_{r+1}'\varep_a(\ei)&= \left\{\begin{array}{ll}(y_r'\psi_r'+1)\varep_a(\ei),&\mbox{if $i_r=i_{r+1}$}\\
y_r'\psi_r'\varep_a(\ei),&\mbox{if $i_r\neq i_{r+1}$}\end{array}\right.\\
y_{r+1}'\psi_r'\varep_a(\ei)&= \left\{\begin{array}{ll}(\psi_r'y_r'+1)\varep_a(\ei),&\mbox{if $i_r=i_{r+1}$}\\\psi_r'y_r'\varep_a(\ei),&\mbox{if $i_r\neq i_{r+1}$}\end{array}\right.\\
(\psi_r')^2\varep_-(\ei)&= \left\{\begin{array}{ll} 0,&\mbox{if $i_r=i_{r+1}$}\\
(y_r'-y_{r+1}')\varep_+(\ei),&\mbox{if $i_r\to i_{r+1}$}\\
(y_{r+1}'-y_r')\varep_+(\ei),&\mbox{if $i_r\leftarrow i_{r+1}$}\\
\varep_-(\ei),&\mbox{otherwise}\end{array}\right.\\
\psi_r'\psi_{r+1}'\psi_r'\varep_+(\ei)&=\left\{\begin{array}{ll} \psi_{r+1}'\psi_r'\psi_{r+1}'\varep_+(\ei)+\varep_-(\ei),&\mbox{if $i_{r+2}=i_r\to i_{r+1}$}\\
\psi_{r+1}'\psi_r'\psi_{r+1}'\varep_+(\ei)-\varep_-(\ei),&\mbox{if $i_{r+2}=i_r\leftarrow i_{r+1}$}\\
\psi_{r+1}'\psi_r'\psi_{r+1}'\varep_+(\ei),&\mbox{otherwise}\end{array}\right.
\end{align*}
for $\ei,\ej\in I^\alpha$, $a,b\in C_2$ and all admissible $r$ and $s$. 
\end{Definition}

%\AM{Change to working with blocks above for finiteness}

\begin{Remark}$\;$
\begin{itemize}
\item[(i)]  As usual, the algebra $\mathcal{R}_{[\alpha]}'$ is a block of the full signed quiver Hecke algebra $\mathcal{R}_n'$, which decomposes as 
\[
\mathcal{R}_n(\Gamma)'=\mathcal{R}_n'(\Gamma,\Zo)=\bigoplus_{[\alpha]\in (Q^+_{\Gamma})_\tau}\mathcal{R}_{[\alpha]}'
\]
and $\mathcal{R}_\alpha'=\varep_{\alpha}^+\mathcal{R}_n'$ where $\varep_{\alpha}^+=\sum_{\ei\in I^\alpha}\varep_+(\ei)$. 
\item[(ii)] The generators in Definition \ref{twisteddef} are somewhat superfluous: for example we can use the relations to see that $\varepsilon_a(\ei)=0$ whenever $\ei=\tau(\ei)$. These extra generators do however mean we can more easily compare the algebra $\mathcal{R}'_{[\alpha]}$ with the quiver Hecke algebra, which we are about to do. 
\end{itemize}

\end{Remark}

The proof of the following lemma requires nothing more than an inspection of the relations in Definition \ref{twisteddef} which we leave to the reader (see \S2.1 for general remarks on algebras graded by a finite group $G$). 

\begin{Lemma}\label{z2grading}
The relations in Definition \ref{twisteddef} are homogeneous with respect to the degree function $\deg_2:\mathcal{R}_n'\to (\Z\times C_2)$ given by 
\[
\deg_2  \varep_a(\ei)= (0,a),\qquad \deg_2 y_r'=(2,-),\qquad \begin{array}{l}\deg_2\psi_r'\varep_+(\ei)= (-c_{i_ri_{r+1}},-)\\\deg_2\psi_r'\varep_-(\ei)= (-c_{i_ri_{r+1}},+)\end{array}
\]
for all $\ei\in I^n$, $a\in C_2$ and $1\leq r\leq n$. In particular, $\mathcal{R}'_n$ is a $(\Z\times C_2)$-graded $\Zo$-algebra.
\end{Lemma}

\begin{Remark}
By forgetting the $C_2$-grading (but {\em not} the $\Z$-grading) on $\mathcal{R}'_n$ we obtain a $\Z$-graded algebra. 
\end{Remark}

\begin{Proposition}
Let $n> 0$ and suppose $e> 2$. If $2$ is invertible in $\Zo$, then for $[\alpha]\in (Q^+_\Gamma)_\tau$, 
\[
\mathcal{R}'_{[\alpha]}(\Gamma,\Zo)\cong \bigoplus_{\beta\in[\alpha]}\mathcal{R}_{\beta}(\Gamma,\Zo)
\]
as $\Z$-graded $\Zo$-algebras. 
\end{Proposition}
\begin{proof}
Define a map $\vartheta:\mathcal{R}'_{[\alpha]}(\Gamma,\Zo)\to \bigoplus_{\beta\in[\alpha]}\mathcal{R}_{\beta}(\Gamma,\Zo)$ on generators by
\begin{align*}
\vartheta(y_r')&= y_r\\
\vartheta(\psi_s')&= \psi_s\\
\vartheta(\varep_a(\ei))&= \left\{\begin{array}{ll}e(\ei)+e(\tau(\ei)),&\mbox{if $a=+$}\\e(\ei)-e(\tau(\ei)),&\mbox{if $a=-$}\end{array}\right.
\end{align*}
for $1\leq r\leq n$ and $1\leq s<n$ and $\ei\in I^\alpha$. We must check $\vartheta$ is an algebra homomorphism of degree zero; this amounts to the largely tedious and straightforward calculation of checking it preserves the relations. We check two relations so the reader can obtain a taste for how they are done. For example, let $\ei,\ej\in I^\alpha$. Then we have
\begin{align*}
\vartheta(\varep_a(\ei)\varep_b(\ej))&= \vartheta(\delta_{\ei\ej}\varep_{ab}(\ei))\\
&= \left\{\begin{array}{ll}\delta_{\ei\ej}(e(\ei)+e(\tau(\ei))),&\mbox{if $a=b$}\\\delta_{\ei\ej}(e(\ei)-e(\tau(\ei))),&\mbox{if $a\neq b$}\end{array}\right.\\
&= \left\{\begin{array}{ll}\delta_{\ei\ej}(e(\ei)+e(\tau(\ei))),&\mbox{if $a=b$}\\\delta_{\ei\ej}(e(\ei)-e(\tau(\ei))),&\mbox{if $a\neq b$}\end{array}\right.\\
&= \vartheta(\varep_a(\ei))\vartheta(\varep_b(\ej)). 
\end{align*}
Now suppose that $\ei\in I^\alpha$ is such that $i_{r+2}=i_r\to i_{r+1}$. Then
\begin{align*}
\vartheta(\psi_r'\psi_{r+1}'\psi_r'\varep_+(\ei))&= \vartheta(\psi_{r+1}'\psi_r'\psi_{r+1}'\varep_+(\ei)+\varep_-(\ei))\\
&= \psi_{r+1}\psi_r\psi_{r+1}[e(\ei)+e(\tau(\ei))]+e(\ei)-e(\tau(\ei))\\
&= (\psi_{r+1}\psi_r\psi_{r+1}+1)e(\ei)+(\psi_{r+1}\psi_r\psi_{r+1}-1)e(\tau(\ei))\\
&= \psi_r\psi_{r+1}\psi_r(e(\ei)+e(\tau(\ei)))\\
&= \vartheta(\psi_r')\vartheta(\psi_{r+1}')\vartheta(\psi_r')\vartheta(\varep_+(\ei))
\end{align*}
as required. We leave the remaining checks to the reader; this amounts to proving $\vartheta$ is surjective. Similarly, define a map $\varsigma:\bigoplus_{\beta\in[\alpha]}\mathcal{R}_\beta\to \mathcal{R}_{[\alpha]}'$ on generators by
\begin{align*}
\varsigma(y_r)&= y_r'\\
\varsigma(\psi_s)&= \psi_s'\\
\varsigma(e(\ei))&= \frac12\left(\varep_+(\ei)+\varep_-(\ei)\right)
%\left\{\begin{array}{ll}\frac12\left(\varep_+(\ei)+\varep_-(\ei)\right),&\mbox{if $\ei\in I^n_+$}\\\frac12\left(\varep_+(\ei)-\varep_-(\ei)\right),&\mbox{if $\ei\in I^n_-$}\end{array}\right.
\end{align*}
for $1\leq r\leq n$, $1\leq s< n$ and $\ei\in\bigcup_{\beta\in [\alpha]}I^\beta$. We must check that $\varsigma$ extends to an algebra homomorphism, again by checking it preserves all relations: note that since $\varep_a(\ei)\varep_b(\ei)=\varep_{a+b}(\ei)$ for all $\ei\in I^n$ and $a,b\in\Z_2$, we can multiply all the relations in Definition \ref{twisteddef} by $\varep_-(\ei)$ to obtain a list of additional relations which also hold in $\mathcal{R}_{[\alpha]}'$, obtaining relations like
\[
\psi_r'y_s'\varep_+(\ei)=y_s'\psi_r'\varep_+(\ei).
\]
These, together with the original list of relations, allow one to check that $\varsigma$ respects all of the relations. For example, if we compute
\[
\varsigma(y_r\psi_se(\ei))=\frac12\Bigl(y_r'\psi_s'\varep_+(\ei)+y_r'\psi_s'\varep_-(\ei)\Bigr)
\] 
we can see that both terms on the right-hand side are indeed relations in $\mathcal{R}'_{[\alpha]}$ (with the first term being a relation from Definition \ref{twisteddef} multiplied on the right by $\varep_-(\ei)$). 
Since $\varsigma$ is also homogeneous of degree zero and since we clearly have $\vartheta\circ\varsigma=\mathrm{id}$ and $\varsigma\circ\vartheta=\mathrm{id}$, this establishes the required isomorphism of $\Z$-graded algebras. 
\end{proof}

%\CB{There should probably be a comment in the proof about how the $\varep_a(\ei)=a\varep_{a}(-\ei)$ gets turned into the equivalence relations.}

%\CB{How does this isomorphism work without mention of $I^n_+$ and $I^n_-$?}

Using the $C_2$-grading, we can write the decomposition of $\mathcal{R}_n'$ into odd and even parts afforded by its $C_2$ grading as $\mathcal{R}_n'=(\mathcal{R}_n')_+\oplus (\mathcal{R}_n')_-$. This, combined with Theorem \ref{aqHbasis}, gives the following corollary regarding the alternating quiver Hecke algebra. 

\begin{Corollary}\label{Z2decomp}
Let $n\geq 0$, let $\Gamma$ be a simply-laced quiver, and let $2$ be invertible in $\Zo$. Then $(\mathcal{R}_n(\Gamma)')_+\cong \mathcal{R}_n(\Gamma/\Gamma')^\sgn$.
%\begin{itemize}
%\item[(i)] If $\alpha\in Q^+_{\Gamma/\Gamma'}$, $(\mathcal{R}_{\alpha}'\oplus\mathcal{R}_{\alpha'}')_+(e,\O)\cong \mathcal{R}_{\alpha/\alpha'}^\sgn(e,\O)$ as $\Z$-graded $\O$-algebras.
%\item[(ii)] If $\alpha\in Q^\pm_{\Gamma/\Gamma'}$, $(\mathcal{R}_\alpha')_+(e,\O)\cong \mathcal{R}_\alpha^\sgn(e,\O)$ as $\Z$-graded $\O$-algebras. 
%\item[(iii)] If $\alpha=n\alpha_0$, $(\mathcal{R}_{n\alpha_0}')_+(e,\O)\cong \mathcal{R}_{n\alpha_0}^\sgn (e,\O)$. 
%\end{itemize} 
%In particular, $(\mathcal{R}_n')_+\cong \mathcal{R}_n^\sgn$. 
\end{Corollary}
%\begin{proof}
%%For parts (i) and (ii), we can equip $\mathcal{R}_n$ with a $C_2$-grading using the Clifford system from Proposition \ref{clifforddecomp}; since the alternating quiver Hecke algebra is the even part of this algebra under this decomposition, the result follows by comparing the even parts of the algebras on both sides. %For part (iii), 
%This follows immediately by comparing the even parts of the algebras on both sides; in one case, the algebra is graded by the $C_2$-grading defined in Lemma \ref{z2grading}, in the other, a $C_2$-grading is inherited from the $\sgn$ involution. The final statement follows by using the direct sum decomposition \eqref{aqHasummands}. 
%\end{proof}

%\CB{Also need to show this for $n\alpha_0$; I think this proof can be redone so it doesn't actually require the Clifford decomposition...}

%\CB{This needs to be changed to be blockwise. Also, importantly, {\em this doesn't work when $\alpha=n\alpha_0$}.}

\renewcommand{\Y}{\mathcal{Y}}

%We define a $C_2$-grading $\deg_2$ on $\mathcal{R}_n$ by defining
%\begin{align*}
%\deg_2e^+(\ei)&= 0\\
%\deg_2e^-(\ei)&=1\\
%\deg_2y_r&= 1\\
%\deg_2\psi_r&= 1. 
%\end{align*}
%
%\begin{proposition}
%The $(\Z\times C_2)$-graded algebra with generators $\{\Psi_r,\mathcal{Y}_r,e^\pm(\ei)\}$ subject to the relations
%\begin{align*}
%e^{a}(\ei)e^{b}(\ej)&=\delta_{\ei\ej}e^{ab}(\ei)\\
%\sum_{\ei\in I^n_+}e^+(\ei)&= 1\\
%e^a(\ei)&=ae^a(-\ei)\\
%\psi_re^a(\ei)&= e^a(s_r\cdot\ei)\psi_r\\
%y_re^a(\ei)&= e^a(\ei)y_r
%\end{align*}
%is isomorphic, when forgetting the $C_2$-grading, to $\mathcal{R}_n$. 
%\end{proposition}
%
%\begin{proposition}
%$\mathcal{R}_n^\sgn$ is isomorphic to the even part of this algebra. 
%\end{proposition}

\begin{Theorem}[Generators and relations for alternating quiver Hecke algebras]\label{aqHapres}
Let $n>1$, let $(\Gamma,\tau)$ be a reversible quiver and suppose 2 is invertible in $\Zo$. 
%\begin{itemize}
%\item[(i)] 
Let $[\alpha]\in Q^+_\tau$. Then the alternating quiver Hecke algebra $\mathcal{R}_\alpha(\Gamma/\Gamma')^\sgn$ is generated by the elements
\[
\{\Psi_1,\Psi_2,\ldots,\Psi_{n-1}\}\cup\{\Y_1,\Y_2,\ldots,\Y_n\}\cup \{e[\ei]\mid [\ei]\in I^{[\alpha]}_\tau\}
\]
subject to the relations
\begin{align*}
e[\ei]e[\ej]&= \delta_{[\ei][\ej]}e[\ei],\qquad\qquad \sum_{[\ei]\in I^{[\alpha]}_\tau}e[\ei]= 1\\
%e^+(\ei)&=-e^+(-\ei)\\
\Y_r e[\ei]&= e[\ei]\Y_r\\
\Psi_r e[\ei]&=e[s_r\cdot\ei]\Psi_r\\%\left\{\begin{array}{ll},&\mbox{if $s_r\cdot\ei\in I^n_+$}\\ -e^+(-s_r\cdot\ei)\Psi_r,&\mbox{if $s_r\cdot\ei\in I^n_-$}\end{array}\right.\\
\Y_r\Y_s&= \Y_s\Y_r\\
\Psi_r\Y_se[\ei]&=\Y_s\Psi_re[\ei]\quad\mbox{if $s\neq r,r+1$}\\ %\left\{\begin{array}{ll}\Y_s\Psi_re^+(\ei),&\mbox{if $s\neq r,r+1$ and $s_r\cdot\ei\in I^n_+$}\\ -\Y_s\Psi_re^+(\ei),&\mbox{if $s\neq r,r+1$ and $s_r\cdot\ei\in I^n_-$}\end{array}\right.\\
\Psi_r\Psi_se[\ei]&= \Psi_s\Psi_re[\ei],\quad\mbox{if $\left|r-s\right|>1$}\\
%\end{align*}%\left\{\begin{array}{ll}\Psi_s\Psi_re^+(\ei),&\mbox{if $\left|r-s\right|>1$ and $s_r\cdot\ei$ and $s_s\cdot\ei\in I^n_+$}\\ -\Psi_s\Psi_re^+(\ei),&\mbox{if $\left|r-s\right|>1$ and either $s_r\cdot\ei$ or $s_s\cdot\ei\in I^n_+$}\end{array}\right.\\
%\begin{align*}
\Psi_r\Y_{r+1}e[\ei]&= \left\{\begin{array}{ll} (\Y_r\Psi_r + 1)e[\ei],&\mbox{if $i_r=i_{r+1}$}\\
\Y_r\Psi_r e[\ei],&\mbox{if $i_r\neq i_{r+1}$}\end{array}\right.\\
\Y_{r+1}\Psi_re[\ei]&= \left\{\begin{array}{ll} (\Psi_r\Y_r+1)e[\ei],&\mbox{if $i_r=i_{r+1}$}\\
\Psi_r\Y_re[\ei],&\mbox{if $i_r\neq i_{r+1}$}\end{array}\right.
\end{align*}
\begin{align*}
\Psi_r^2e[\ei]&= \left\{\begin{array}{ll} %\Y_3e^+(\ei),&\mbox{if $s_r\cdot\ei\in I^n_-$ and $i_r\to i_{r+1}$}\\
%-e^+(\ei),&\mbox{if $s_r\cdot\ei\in I^n_-$ and $i_r\leftarrow i_{r+1}$}\\
0,&\mbox{if $i_r=i_{r+1}$}\\
(\Y_r-\Y_{r+1})e[\ei],&\mbox{if $i_r\to i_{r+1}$}\\
(\Y_{r+1}-\Y_r)e[\ei],&\mbox{if $i_r\leftarrow i_{r+1}$}\\
e[\ei],&\mbox{otherwise}\end{array}\right.\\
%\end{align*}
%\begin{align*}
\Psi_r\Psi_{r+1}\Psi_re[\ei]&= \left\{\begin{array}{l} (\Psi_{r+1}\Psi_r\Psi_{r+1}-1)e[\ei],\\
\qquad\qquad\qquad\quad\mbox{if $i_r=i_{r+2}\to i_{r+1}$}\\
(\Psi_{r+1}\Psi_r\Psi_{r+1}+1)e[\ei],\\
\qquad\qquad\qquad\qquad\mbox{if $i_r=i_{r+2}\leftarrow i_{r+1}$}\\
\Psi_{r+1}\Psi_r\Psi_{r+1}e[\ei],\\
\qquad\qquad\qquad\quad\mbox{otherwise}\end{array}\right.
\end{align*}
for all $[\ei]\in I^{[\alpha]}_\tau$ and all admissible $r$ and $s$.
%, together with the additional relation
%\[
%\sum_{\ei\in I^{\alpha/\alpha'_\sim}}e^+([\ei])=1.
%\]
Moreover, $\mathcal{R}_\alpha(\Gamma/\Gamma')^\sgn$ is $\Z$-graded with degree function 
\[
\deg\Psi_re[\ei]=-a_{i_r i_{r+1}},\quad \deg\Y_s=2,\quad \deg e[\ei]=0. 
\]
for all $1\leq r<n$, $1\leq s\leq n$ and $[\ei]\in I^{[\alpha]}_\tau$. 
%\item[(ii)] Let $\alpha\in Q^\pm_{\Gamma/\Gamma'}$. Then the alternating quiver Hecke algebra $\mathcal{R}_\alpha^\sgn$ is generated by
%\[
%\{\Psi_1,\Psi_2,\ldots,\Psi_{n-1}\}\cup\{\Y_1,\Y_2,\ldots,\Y_n\}\cup \{e^+([\ei])\mid [\ei]\in I^{\alpha}_\sim\}
%\]
%subject to the relations \eqref{commonrels} for all $[\ei]\in I^\alpha_\sim$ and all admissible $r$ and $s$, together with the additional relation 
%\[
%\sum_{\ei\in I^\alpha_\sim}e^+([\ei])=1.
%\]
%Moreover, $\mathcal{R}_\alpha^\sgn$ is $\Z$-graded with degree function
%\[
%\deg\Psi_re^+([\ei])=-a_{i_r i_{r+1}},\quad \deg\Y_s=2,\quad \deg e^+([\ei])=0. 
%\]
%for all $1\leq r<n$, $1\leq s\leq n$ and $[\ei]\in I^{\alpha}_\sim$. 
%\item[(iii)] Let $\alpha=n\alpha_0$. Then the alternating quiver Hecke algebra $\mathcal{R}_{n\alpha_0}^\sgn$ is generated by 
%\[
%
%\]
%\item[(ii)] If $\alpha=n\alpha_0$, the alternating quiver Hecke algebra $\mathcal{R}_{n\alpha_0}^\sgn$ is generated by the elements
%\begin{align*}
%\{y_ry_s\mid 1\leq r<s\leq n\}&\cup\{\psi_ry_s,y_s\psi_r\mid 1\leq r<n, \, 1\leq s\leq n\}\\
%&\qquad\qquad\qquad\;\;\cup\{\psi_r\psi_s\mid 1\leq r,s<n\}\cup\{e\}
%\end{align*}
%subject to the relations
%\begin{align*}
%e&= 1\\
%y_ry_se&= y_sy_re\\
%\psi_ry_se&= y_s\psi_re,\quad\mbox{if $s\neq r,r+1$}\\
%\psi_r\psi_se&= \psi_s\psi_re,\quad\mbox{if $\left|r-s\right|>1$}\\
%\psi_ry_{r+1}e&= y_r\psi_re+e\\
%y_{r+1}\psi_re&= \psi_ry_re+e\\
%\psi_r^2e&= e
%\end{align*}
%for all admissible $r$ and $s$, where $e=\underbrace{(0,0,\ldots,0)}_{\text{$n$ zeroes}}$. 
%%\end{proof}
%\end{itemize}
\end{Theorem}

%\CB{Fix up the equivalence classes so that this can be done blockwise}

\begin{proof}
%By the decomposition in Corollary \ref{Z2decomp}, the alternating quiver Hecke algebra will consist of ``even'' relations from the full $C_2$-twisted quiver Hecke algebra. An inspection of the list of relations demonstrates that the even ones are precisely those included in the above list. 
%We start by proving part (i),
First we note that all the requirements on the residue sequences $\ei$ in fact depend only on the {\em equivalence class} of the sequence because of the symmetry of the Cartan matrix in \eqref{cartanmatrixdef}. By Corollary \ref{Z2decomp}, it is enough to prove that the abstract algebra $A_{[\alpha]}$ defined in the statement of the theorem is isomorphic to $(\mathcal{R}_{[\alpha]}')_+$. Define a map $\varrho:A_{[\alpha]}\to (\mathcal{R}_{[\alpha]}')_+$ by
\[
e[\ei]\mapsto \varepsilon_+(\ei),\quad \mathcal{Y}_r\mapsto y_r'\varepsilon_-(\ei),\quad \Psi_s\mapsto \psi_s'\varepsilon_-(\ei),
\]
for all $[\ei]\in I^{[\alpha]}_\tau$, $1\leq r\leq n$ and $1\leq s<n$. It is a straightforward check that all the relations of $A_{[\alpha]}$ are satisfied in $(\mathcal{R}_{[\alpha]}')_+$, so $\varrho$ determines a well-defined homogeneous algebra homomorphism of degree zero. 

By definition, $(\mathcal{R}_{[\alpha]}')_+$ is generated by even words in the generators of $\mathcal{R}_{[\alpha]}'$. However the only even {\em generators} of $\mathcal{R}_{[\alpha]}'$ are the idempotents $\varepsilon_+(\ei)$ for $\ei\in I^\alpha\cup I^{\alpha'}$, so $(\mathcal{R}_{[\alpha]}')_+$ is generated by these idempotents together with all words of even length in the {\em odd} generators of $\mathcal{R}_{[\alpha]}'$. It is now easy to see that $(\mathcal{R}_{[\alpha]}')_+$ is generated by the images of the generators of $A_{[\alpha]}$ under $\varrho$ and so $\varrho$ is surjective. 

The algebra $\mathcal{R}_{[\alpha]}'$ is defined by generators and relations, so $(\mathcal{R}_{[\alpha]}')_+$ is the subalgebra of $\mathcal{R}_{[\alpha]}'$ generated by the even words in the generators of $\mathcal{R}_{[\alpha]}'$ modulo the even part of the relation ideal defining $\mathcal{R}_{[\alpha]}'$, which is the set of all linear combinations of arbitrary products of even relations multiplied by even products of odd relations. %The even relations in $\mathcal{R}_{[\alpha]}'$ are generated by all even words in the relations in $\mathcal{R}_{[\alpha]}'$. Such a relation will be a product of even relations, and even numbers of odd relations. 
However the only even relations in $\mathcal{R}_{[\alpha]}'$ are given by idempotents and commutation relations.
 %, which are clearly mapped to by the corresponding idempotent and commutation relations in $A_{[\alpha]}$ under $\varrho$. 
Therefore, the even part of the relation ideal for $\mathcal{R}_{[\alpha]}'$ is generated by the even relations in $\mathcal{R}_{[\alpha]}'$ together with all products of the odd relations in $\mathcal{R}_{[\alpha]}'$. 
It follows that all the remaining relations are generated by even products of odd relations, together with odd relations multiplied by $\varepsilon$: in this way we obtain the complete set of relations for $(\mathcal{R}_{[\alpha]}')_+$. 
One checks these are precisely the relations written above for $A_{[\alpha]}$; for example multiplying the relation 
\[
y_{r+1}'\psi_r'\varepsilon_-(\ei)=\left\{\begin{array}{ll}(\psi_r'y_r'+1)\varepsilon_-(\ei),&\mbox{if $i_r=i_{r+1}$}\\ \psi_r'y_r'\varepsilon_-(\ei),&\mbox{if $i_r\neq i_{r+1}$}\end{array}\right.
\]
by $\varepsilon_-(\ei)$ and using the idempotent relations in $\mathcal{R}_{[\alpha]}'$ gives the relation
\[
y_{r+1}'\varepsilon_-(\ei)\psi_r'\varepsilon_-(\ei)\varepsilon_+(\ei)=\left\{\begin{array}{ll}(\psi_r'\varepsilon_-(\ei)y_{r+1}'\varepsilon_-(\ei)+1)\varepsilon_+(\ei),&\mbox{if $i_r=i_{r+1}$}\\ \psi_r'\varepsilon_-(\ei) y_r'\varepsilon_-(\ei)\varepsilon_+(\ei),&\mbox{if $i_r\neq i_{r+1}$}\end{array}\right.
\]
which is precisely the image of the relation
\[
\mathcal{Y}_{r+1}\Psi_re[\ei]=\left\{\begin{array}{ll}(\Psi_r\mathcal{Y}_r+1)e[\ei],&\mbox{if $i_r=i_{r+1}$}\\ \Psi_r\mathcal{Y}_re[\ei],&\mbox{if $i_r\neq i_{r+1}$}\end{array}\right.
\] 
under $\varrho$. Continuing in this way we see $\varrho$ is an isomorphism. 

Finally, we note that the given degree function is well-defined; since the Cartan matrix $(c_{ij})_{i,j\in I}$ is symmetric, the entries only depend on equivalence classes of residue sequences under $\tau$. 

%For part (ii), the given relations are precisely the even relations from Definition \ref{twisteddef} in the particular case when $\alpha=n\alpha_0$ (for example, $i_r=i_{r+1}$ is always true in this case, which gives the simplified quadratic relations), together with the even relations from Definition \ref{twisteddef} obtained by multiplying the odd relations by $\varepsilon$. 
\end{proof}
%\CB{Finish this (make sure the list of relations is right since I've changed notation) -- I think the quadratic relation is too complicated, and maybe the mixed $\Psi$ and $\mathcal{Y}$ relations as well?}

%\AM{You also need the missing ``even'' relaitons implied in 5.2.19 given by multiplying by $\varepsilon$.}

\begin{Remark}$\;$

\begin{itemize}
\item[(i)] Since we have given an isomorphism between the alternating quiver Hecke algebra and the even part of the abstract algebra $\mathcal{R}_n'$, for which we gave an abstract presentation by generators and relations, we see that the generators from Definition \ref{aqHagens} do not depend on the choice of equivalence class representatives $\ei_\tau$ (see p\pageref{eiplusdef}). 
\item[(ii)] In \cite[\S5.4]{BoysThesis} and \cite{BM}, cyclotomic quotients of the algebras from this paper are studied. In particular, \cite[Theorem B]{BM} gives a presentation for a cyclotomic quotient of the alternating quiver Hecke algebras by generators and relations very similar to those in this paper, in the special case when $\Gamma=\Gamma_e$ and the cyclotomic relation is simple. In the cyclotomic case when the cyclotomic parameter $\Lambda=\Lambda_0$, it is shown in \cite[Chapter 6]{BoysThesis} that alternating cyclotomic quiver Hecke algebras are isomorphic to alternating Hecke algebras in the sense of Mitsuhashi \cite{Mitsuhashi:A}, giving an analogue of Brundan and Kleshchev's result \cite{BK:GradedKL}. 
\end{itemize}
\end{Remark}

    \let\origHref=\href
    \renewcommand\href[3][\relax]{#3}

\bibliographystyle{andrew}
%\bibliography{papers}

    \let\href=\origHref

\end{document}